\font\tenfrakturb=eufb10
\font\tenfraktur=eufm10
\font\tenmsbm=msbm10
\font\sevenfrakturb=eufb7
\font\sevenfraktur=eufm7
\font\sevenmsbm=msbm7
\font\fivefrakturb=eufb5
\font\fivefraktur=eufm5
\font\fivemsbm=msbm5
\newfam\bgothicfam
\newfam\gothicfam
\newfam\msbmfam
\textfont\bgothicfam = \tenfrakturb \scriptfont\bgothicfam=\sevenfrakturb
\scriptscriptfont\bgothicfam=\fivefrakturb

\textfont\gothicfam = \tenfraktur \scriptfont\gothicfam=\sevenfraktur
\scriptscriptfont\gothicfam=\fivefraktur

\textfont\msbmfam = \tenmsbm \scriptfont\msbmfam=\sevenmsbm
\scriptscriptfont\msbmfam=\fivemsbm

\def\frak{\tenfraktur\fam\gothicfam}
\def\Bbb{\tenmsbm\fam\msbmfam}

\catcode`@=11
\def\renewcounter#1{\@definecounter{#1}\@ifnextchar[{\@newctr{#1}}{}}

\documentstyle[12pt]{article}
\renewcounter{equation}[section]

\newtheorem{th}{Theorem}[section]

\begin{document}
\title{Asymptotics of the Heat Kernel on Rank 1 Locally Symmetric Spaces}
\author{A.A. Bytsenko \thanks{E-mail: abyts@fisica.uel.br\,\,\, On leave from
Sankt-Petersburg State Technical University}\\
Departamento de Fisica, Universidade Estadual de Londrina,\\
Caixa Postal 6001, Londrina-Parana, Brazil\\
and \\
\ F.L. Williams
\thanks{E-mail: williams@math.umass.edu}\\
Department of Mathematics, University of Massachusetts,\\
Amherst, Massachusetts 01003}

\date{April, 1998}
\maketitle
\begin{abstract}

We consider the heat kernel (and the zeta function) associated with 
Laplace type operators acting on a general irreducible rank 1 locally
symmetric space $X$. The set of Minakshisundaram-Pleijel coefficients 
$\{A_k(X)\}_{k=0}^{\infty}$ in the short-time asymptotic expansion of the 
kernel is calculated explicitly.

\end{abstract}

\section{Introduction}

In the theory of quantum fields on curved background spaces, the short-time
expansion of the heat kernel plays an extremely important role. In particular
situations, for example, the coefficients in the expansion control the
one-loop divergences of the effective action, and related quantities such
as the stress energy momentum tensor. Some of these coefficients have been
determined and appear in the physics and mathematical literature. Note
the references \cite{miat76,miat79-29-249,miat80-260-1,dewi65,seel67-10-288,
atiy73-19-279} for closed Riemannian
manifolds and \cite{moss89-229-261,bran90-15-245} for manifolds with a smooth
boundary. The literature on these matters is very vast.

In Refs.  \cite{miat76,miat79-29-249,miat80-260-1}, R. Miatello studies the 
case of a closed locally
symmetric rank 1 manifold $X$, using the representation theory of the group
of isometries of $X$. We consider the same case in the present paper, but
we use the spectral zeta function of $X$. By our approach we determine the
expansion coefficients immediately and explicitly (essentially in one step),
given the results of \cite{will98-182-137}.
Recently the topological Casimir energy \cite{will97-38-796}, the one-loop
effective action, and the multiplicative and conformal anomaly
\cite{byts97u-20,byts98-13-99} associated with Laplace type operators on $X$,
and their product, have been analysed also by use of the spectral zeta
function.

The paper is organized as follows. In Sect. 2 we define the spectral zeta
function $\zeta_{\Gamma}(s;\chi)$ of $X$ corresponding to a finite-dimensional
representation $\chi$ of the fundamental group $\Gamma$ of $X$. The residues 
of $\zeta_{\Gamma}(s;\chi)$ and special values of this zeta function, which
relate to the expansion coefficients, are provided by Theorems 2.1 and 2.2
In Sect. 3 we consider the asymptotic expansion of the heat kernel (as 
$t\rightarrow 0^{+}$), and compute all the expansion coefficients in closed 
form in the main theorem, Theorem 3.1. Sect. 4 contains some remarks in 
summary. We include an Appendix with information supplementary to 
Theorems 2.1, 2.2 and 3.1.

\section{The Spectral Zeta Function}

We shall be working with an irreducible rank 1 symmetric space $M=G/K$
of non-compact type. Thus $G$ will be a connected non-compact simple split rank
1 Lie group with finite center and $K\subset G$ will be a maximal compact
subgroup \cite{helg62}. Let $\Gamma\subset G$ be a discrete,
co-compact torsion free subgroup. Then $X=X_{\Gamma}=\Gamma\backslash M$ is a 
compact Riemannian manifold with fundamental group $\Gamma$; namely $X$ is a 
compact locally symmetric space. Given a finite-dimensional unitary 
representation $\chi$ of $\Gamma$ there is the corresponding vector bundle 
$V_{\chi}\mapsto X$ over $X$ given by $V_{\chi}=\Gamma\backslash (M\otimes 
F_{\chi})$, where $F_{\chi}$ (the fibre of $V_{\chi}$) is the representation 
space of $\chi$ and where $\Gamma$ acts on $M\otimes F_{\chi}$ by the rule
$\gamma\cdot(m,f)=(\gamma\cdot m,\chi(\gamma)f)$ for $(\gamma,m,f)\in (\Gamma
\otimes M\otimes F_{\chi})$. Let $\Delta_{\Gamma}$ be the Laplace-Beltrami
operator of $X$ acting on smooth sections of $V_{\chi}$; we obtain 
$\Delta_{\Gamma}$ by projecting the Laplace-Beltrami operator of $M$ (which is
$G-$ invariant and thus $\Gamma-$ invariant) to $X$. As $X$ is compact we can
consider the spectrum $\{\lambda_j=\lambda_j(\chi),\,\,
n_j=n_j(\chi)\}_{j=0}^{\infty}$ of $-\Delta_{\Gamma}$, where $n_j$ is
the (finite) multiplicity of the eigenvalue $\lambda_j$. We use the minus 
preceding $\Delta_{\Gamma}$ to have the $\lambda_j\geq 0:\,\,0=\lambda_0
<\lambda_1<\lambda_2\,\, ...;$\,\, $\lim_{j\rightarrow\infty}\lambda_j=\infty$.

The spectral zeta function $\zeta_{\Gamma}(s;\chi)$ of $X_{\Gamma}$ of
Minakshisundaram-Pleijel type \cite{mina49-1-242}, which we shall consider is 
defined by

\begin{equation}
\zeta_{\Gamma}(s;{\chi})=\sum_{j=1}^{\infty}\frac{n_j(\chi)}
{\lambda_j(\chi)^{s}}\,
\mbox{,}
\end{equation}
for $\Re s>>0$. $\zeta_{\Gamma}(s;\chi)$ is a holomorphic function on the
domain $\Re s>d/2$, where $d$ is the dimension of $M$, and by general 
principles $\zeta_{\Gamma}(s;\chi)$ admits a meromorphic continuation to the
full complex plane ${\Bbb C}$. However since the manifold $X_{\Gamma}$ is
quite special it is desirable to have the meromorphic continuation of
$\zeta_{\Gamma}(s;\chi)$ in an explicit form, for example in terms of the
structure of $G$ and $\Gamma$. Using the Selberg trace formula and the $K$-
spherical harmonic analysis of $G$, we have obtained such a form in 
\cite{will98-182-137}; also see Refs. \cite{rand75-201-241,will97-38-796}. 
In particular we can obtain the residues of
$\zeta_{\Gamma}(s;\chi)$, and compute the special values 
$\zeta_{\Gamma}(-n;\chi)$, \,\,$n=0, 1, 2, ... -$ results which play a
decisive role in the present work. To state these results we introduce further
notation.

Up to local isomorphism we can represent $M=G/K$ by the following quotients:

\begin{equation}
M=\left[ \begin{array}{ll}SO_1(n,1)/SO(n)\,\,\,\,\,\,\,\,\,\,\,\,\,\,\,\,\,
\,\,\,\,\,\,\,\,\,\,\,\,\,\,\,\,\,\,\,\,\,(I) \\
SU(n,1)/U(n)\,\,\,\,\,\,\,\,\,\,\,\,\,\,\,\,\,\,\,\,\,\,\,\,\,\,\,\,\,\,\,\,
\,\,\,\,\,\,\,\,\,(II)
\\SP(n,1)/(SP(n)\otimes SP(1))\,\,\,\,\,(III)\\
F_{4(-20)}/Spin(9)\,\,\,\,\,\,\,\,\,\,\,\,\,\,\,\,\,\,\,\,\,\,\,\,\,\,\,\,\,
\,\,\,\,\,\,(IV)
\end{array} \right]
\mbox{,}
\end{equation}
where $d=n, 2n, 4n, 16$ respectively. We shall need the real number $\rho_0$
which corresponds to $1/2$ the sum of the positive real restricted roots of
$G$ with respect to a nilpotent factor in an Iwasawa decomposition of $G$.
$\rho_0$ is given by $\rho_0=(n-1)/2, n, 2n+1, 11$ respectively in the cases
$(I)$ to $(IV)$. For details on these matters the reader may consult 
\cite{helg62}, and also the Appendix in \cite{will97-38-796}.

The spherical harmonic analysis on $M$ is controlled by Harish-Chandra's 
Plancherel density $|C(r)|^{-2}$, a function on the real numbers $\Bbb R$,
computed by Miatello \cite{miat76,miat79-29-249,miat80-260-1}, and others, in 
the rank 1 case we are
considering. We choose a normalization of Haar measure on $G$ however which
differs from that of \cite{miat76,miat79-29-249,miat80-260-1}; see Ref. 
\cite{will98-182-137}. For a suitable constant $C_G$
depending only on $G$, and for a suitable even polynomial $P(r)$ of degree
$d-2$ for $G\neq SO_1(n,1)$ with $n$ odd, and of degree $d-1=2m$ for 
$G=SO_1(2m+1,1)$, \,\,$|C(r)|^{-2}$ is given by 

\begin{equation}
|C(r)|^{-2}=\left[ \begin{array}{ll}C_G\pi rP(r)\tanh(\pi r) \hspace{1.0cm}
\mbox{for $G=SO_1(2m,1)$},\\
C_G\pi rP(r)\tanh(\pi r/2) \hspace{0.6cm}\mbox{for $G=SU
(n,1),\,\,\,\,\,\,\, n$ odd}\\
\hspace{4,5cm}\mbox{or $G=SP(n,1),\,\,\,\,\, F_{4(-20)}$},\\C_G\pi
rP(r)\coth(\pi r/2) \hspace{0.7cm}\mbox{for $G=SU
(n,1),\,\,\,\,\,\,n$ even},\\
C_G\pi P(r) \hspace{2.8cm}\mbox{for $G=SO_1(2m+1,1)$}\\
\end{array} \right]
\mbox{.}
\end{equation}
The value of $C_G$ and the explicit form of $P(r)$ is given in the Appendix.
For real hyperbolic space $M=SO_1(2m,1)/SO(2m)$ of even dimension $2m$,
for example, $P(r)$ is given by

\begin{equation}
P(r)=\prod_{j=0}^{m-2}\left[r^2+\frac{(2j+1)^2}{4}\right]
\mbox{.}
\end{equation}
The coefficients of $P(r)$ will be denoted by $a_{2j}$:

\begin{eqnarray}
P(r) &=& \sum_{j=0}^{d/2-1}a_{2j}r^{2j} \hspace{1.0cm}
\mbox{for}\,\,\, G\neq SO_1(2m+1,1), \nonumber\\
&=& \sum_{j=0}^m a_{2j}r^{2j} \hspace{1.3cm}\mbox{for}\,\,\, G=SO_1(2m+1,1).
\end{eqnarray} 
 We denote by ${\rm Vol}(\Gamma\backslash G)$ the $G-$ invariant volume of 
$\Gamma\backslash G$ induced by Haar measure on $G$.

As pointed out earlier, the explicit meromorphic structure of the zeta
function $\zeta_{\Gamma}(s;\chi)$ of (2.1) is worked out in 
\cite{will98-182-137}
in terms of the spherical harmonic analysis of $G$ and $\Gamma-$ structure;
see Theorems 4.2, 5.1 there; also compare Theorems 5.2, equation (6.1), and
Theorem 6.9 of \cite{will97-38-796}. In particular, apart from the case 
$G=SO_1(n,1)$
with $n$ odd (a case which we treat separately), $\zeta_{\Gamma}(s;\chi)$
is holomorphic except for possibly simple poles at $s=1, 2, ..., d/2$.
By Theorem 5.1 of \cite{will98-182-137}, or by the results stated in 
\cite{will97-38-796} we can compute the residues at these points 
$s=1, 2, ..., d/2$. The results are the following, where we omit the 
cotangent case.

\begin{th}
Apart from the cases $SO_1(\ell,1)$,\,\,  $SU(q,1)$ with $\ell$ odd and $q$
even, the residue of $\zeta_{\Gamma}(s ;\chi)$ at $s=m$ (for $m$ an integer,
$1\leq m\leq d/2$) equals

\begin{equation}
\frac{1}{4}\chi(1)\mbox{Vol}(\Gamma\backslash G)
C_G\sum_{j=0}^{d/2-m}(-1)^j
\left( \begin{array}{ll}
m+j-1\\
\,\,\,\,\,\,\,\,\,\,\,\,\, j \\
\end{array} \right)
\rho_0^{2j}a_{2(m+j-1)}
\mbox{,}
\end{equation}
given the preceding notation. Also for $n=1, 2, ...,$

\begin{eqnarray}
\zeta_{\Gamma}(-n;\chi) &=& \frac{1}{4}\chi(1)\mbox{Vol}(\Gamma\backslash G)
C_G\left[\sum_{j=0}^{d/2-1}\frac{(-1)^{j+1}j!\rho_0^{2(j+n+1)}a_{2j}}
{(n+1)(n+2)\cdot\cdot\cdot (n+j+1)}\right. \nonumber\\
&+& \left. 2\sum_{j=0}^{d/2-1}\sum_{k=0}^n\frac{(-1)^kn!}{(n-k)!}
\rho_0^{2(n-k)}b_{k+1}(j)a_{2j}\right]
\mbox{,}
\end{eqnarray}
where 

\begin{equation}
b_p(j)\stackrel{def}{=}\left[2^{1-2(p+j)}-1\right]\left[\frac{\pi}{a(G)}
\right]^{2(p+j)}
\frac{(-1)^jB_{2(p+j)}}{2(p+j)[(p-1)!]}
\mbox{,}
\end{equation}
for $p=1, 2, ...,$ $B_r$ the $r$-th Bernoulli number, and for

\begin{equation}
a(G)\stackrel{def}{=}\left[ \begin{array}{ll}\pi \hspace{1.0cm}
\mbox{if\,\, $G=SO_1(\ell,1)$ with $\ell$ even},\\
\frac{\pi}{2} \hspace{1.0cm}\mbox{if\,\, $G=SU(q,1)$\,\,\, with $q$ odd}\\
\hspace{1,2cm}\mbox{or $G=SP(\ell,1)$, any $\ell$, $F_{4(-20)}$}\\
\end{array} \right]
\mbox{.}
\end{equation}
$\zeta_{\Gamma}(0;\chi)=-n_0(\chi)+$ (the right hand side of Eq. (2.7) 
evaluated at $n=0$).
\end{th}            

Now we consider the case $G=SO_1(\ell,1)$ with $\ell$ odd. By the results of 
\cite{will98-182-137}, for $G=SO_1(2n+1,1)$ \,$\zeta_{\Gamma}(s;\chi)$ has at 
most a simple pole at the points $s=d/2-k$, $k=0, 1, 2, ...$. Moreover

\begin{th}
For $G=SO_1(2n+1,1)$ the residue of $\zeta_{\Gamma}(s;\chi)$ at $s=d/2-k$
(where $d/2=n+1/2,\,\,  k=0, 1, 2,...$) equals

\begin{equation}
\frac{1}{4}\chi(1)\mbox{Vol}(\Gamma\backslash G)
C_G\sum_{j=0}^n\frac{(-1)^{j+n+k}\rho_0^{2(j+k-n)}\Gamma(j+\frac{1}{2})a_{2j}}
{(j-n+k)!\Gamma(n+\frac{1}{2}-k)}
\mbox{,}
\end{equation}
for $k\geq n$, and equals

\begin{equation}
\frac{1}{4}\chi(1)\mbox{Vol}(\Gamma\backslash G)
C_G\sum_{j=0}^k\frac{(-1)^j\rho_0^{2j}\Gamma(n-k+j+\frac{1}{2})a_{2(n-k+j)}}
{j!\Gamma(n+\frac{1}{2}-k)}
\mbox{,}
\end{equation}
for $0\leq k <n$. Here $\rho_0=n$. Also \,\, $\zeta_{\Gamma}(0;\chi)=
-n_0(\chi)$, whereas $\zeta_{\Gamma}(-k;\chi)=0$ for $k=1, 2, ...$.
\end{th}

In Theorems 2.1 and 2.2 the constant $C_G$ is given in the Appendix.

\section{The Heat Kernel Coefficients}

The object of interest is the heat kernel $\omega_{\Gamma}(t;\chi)$
defined for $t>0$ by

\begin{equation}
\omega_{\Gamma}(t;\chi)=\sum_{j=0}^{\infty}n_j(\chi)
e^{-\lambda_j(\chi)t}
\mbox{.}
\end{equation}
If $h_t$ is the fundamental solution of the heat equation on $M$, then
$h_t$ and $\omega_{\Gamma}(t;\chi)$ are related by the Selberg trace
formula (cf. \cite{will98-182-137})

\begin{equation}
\omega_{\Gamma}(t;\chi)=\chi(1)\mbox{Vol}(\Gamma\backslash G)h_t(1)
+\theta_{\Gamma}(t;\chi)
\mbox{,}
\end{equation}
where the theta function $\theta_{\Gamma}(t;\chi)$ is given by Eq. (4.18) of
\cite{will98-182-137} (for $b=0$ there) and where

\begin{equation}
h_t(1)=\frac{1}{4\pi}e^{-\rho_0^2t}\int_{{\Bbb R}}
e^{-r^2t}|C(r)|^{-2}dr
\mbox{.}
\end{equation}

We shall {\em not} need the result (3.2). Our goal is to compute explicitly
all of the coefficients $A_k=A_k(\Gamma,\chi)$ in the asymptotic expansion

\begin{equation}
\omega_{\Gamma}(t;\chi)\simeq (4\pi t)^{-d/2}\sum_{k=0}^{\infty}A_kt^k,
\,\,\,\,\,\mbox {as} \,\, t\rightarrow 0^{+}
\mbox{.}
\end{equation}

Now $\zeta_{\Gamma}(s;\chi)$ and $\omega_{\Gamma}(t;\chi)$ are related by the
Mellin transform:

\begin{equation}
\zeta_{\Gamma}(s;\chi)=\frac{{\frak M}[\omega_{\Gamma}](s)}{\Gamma(s)}=
\frac{1}{\Gamma(s)}\int_0^{\infty}\omega_{\Gamma}(t;\chi)t^{s-1}dt,\,\,\,\,\,
\mbox{for}\,\,\,\, \Re s>\frac{d}{2}
\mbox{.}
\end{equation}
Moreover one knows by abstract generalities 
(cf. \cite{mina49-1-242,voro87-110-439} for example)
that the coefficients $A_k$ are related to residues and special values of
$\zeta_{\Gamma}(s;\chi)$. Namely for $m$ an integer with $1\leq m\leq d/2$,
for d even

\begin{equation}
A_{\frac{d}{2}-m}=(4\pi)^{d/2}\Gamma(m)
\times \left[{\rm residue}\,\,\, {\rm of}\,\,\, \zeta_{\Gamma}(s;\chi)\,\,\,
{\rm at}\,\,\, s=m\right]
\mbox{.}
\end{equation}
Also for a positive integer $n$

\begin{equation}
A_{\frac{d}{2}+n}=\frac{(-1)^n(4\pi)^{d/2}}{n!}\zeta_{\Gamma}(-n;\chi)
\mbox{,}
\end{equation}
whereas

\begin{equation}
A_{\frac{d}{2}}=(4\pi)^{d/2}\left[n_0(\chi)+\zeta_{\Gamma}(0;\chi)\right]
\mbox{.}
\end{equation}
For $G=SO_1(2n+1,1)$ (the only case in which $d$ is odd) we have for 
$k=0, 1, 2, ...$

\begin{equation}
A_k=(4\pi)^{d/2}\Gamma\left(\frac{d}{2}-k\right)
\times \left[{\rm residue}\,\,\, {\rm of}\,\,\, \zeta_{\Gamma}(s;\chi)\,\,\,
{\rm at}\,\,\, s=\frac{d}{2}-k\right]
\mbox{;}
\end{equation}
$d/2=n+1/2$. Hence by Eqs. (3.6) - (3.9), and Theorems 2.1, 2.2 we obtain
the following main result.

\begin{th}
The heat kernel $\omega_{\Gamma}(t;\chi)$ in (3.1) admits an asymptotic
expansion (3.4). More precisely, given any non-negative integer $N$ one has

\begin{equation}
\lim_{t\rightarrow 0^{+}}\left[(4\pi t)^{d/2}\omega_{\Gamma}(t;\chi)-
\sum_{k=0}^NA_k(\Gamma,\chi)t^k\right]t^{-N}=0
\end{equation}
where, apart from the cotangent case in (2.3) (i.e. the case $G=SU(q,1)$ with
q even), the coefficients $A_k(\Gamma,\chi)=A_k(X_{\Gamma})$ are given as 
follows. 

For all $G$ except $G=SO_1(\ell,1)$, $SU(q,1)$ with $\ell$ odd and $q$ even

\begin{eqnarray}
A_k(\Gamma,\chi) &=& (4\pi)^{\frac{d}{2}-1}\chi(1){\rm Vol}
(\Gamma\backslash G)C_G\pi\sum_{\ell=0}^k
\frac{(-\rho_0^2)^{k-\ell}}{(k-\ell)!}[\frac{d}{2}-(\ell+1)]!
a_{2[\frac{d}{2}-(\ell+1)]}
\nonumber\\
&& {\rm for}\,\,\, 0\leq k\leq \frac{d}{2}-1 
\mbox{,}
\end{eqnarray} 

\begin{eqnarray}
A_{\frac{d}{2}+n}(\Gamma,\chi) &=& (-1)^n(4\pi)^{\frac{d}{2}-1}\chi(1)
{\rm Vol}(\Gamma\backslash G)C_G\pi
\left[\sum_{j=0}^{\frac{d}{2}-1}(-1)^{j+1}\frac{\rho_0^{2(n+1+j)}j!a_{2j}}
{(n+1+j)!} \right. \nonumber\\
&+& \left. 2\sum_{j=0}^{\frac{d}{2}-1}\sum_{\ell=0}^n(-1)^{\ell}
\frac{\rho_0^{2(n-\ell)}}{(n-\ell)!}b_{\ell+1}(j)a_{2j}\right] \nonumber\\
&& {\rm for}\,\,\,n=0, 1, 2, ...,
\end{eqnarray}
where $b_p(j)\,\,(p=1,2, ...)$ and $a(G)$ are given by (2.8) and (2.9). 

For $G=SO_1(2n+1,1),\,\, k=0,1,2, ...$
   
\begin{equation}
A_k(\Gamma,\chi)=\pi(4\pi)^{n-\frac{1}{2}}\chi(1){\rm Vol}
(\Gamma\backslash G)C_G
\sum_{\ell=0}^{{\rm min}(k,n)}
\frac{(-n^2)^{k-\ell}\Gamma\left(n-\ell+\frac{1}{2}\right)a_{2(n-\ell)}}
{(k-\ell)!}
\mbox{,}
\end{equation}
or

\begin{equation}
A_k(\Gamma,\chi)=\pi^{3/2}(4\pi)^{n-\frac{1}{2}}\chi(1){\rm Vol}
(\Gamma\backslash G)C_G
\sum_{\ell=0}^{{\rm min}(k,n)}
\frac{(-\rho_0^2)^{k-\ell}[2(n-\ell)]!a_{2(n-\ell)}}{(k-\ell)!(n-l)!
2^{2(n-\ell)}}
\mbox{,}
\end{equation}
using that $\Gamma(m+1/2)=\pi^{1/2}(2m)!\left[2^{2m}m!\right]^{-1}$.
\end{th}

\section{Conclusions}

Using results \cite{will98-182-137} on the meromorphic structure of the zeta 
function of
a rank 1 locally symmetric space $X$, we have obtained in a quick
computation all of the Minakshisundaram-Pleijel coefficients (in closed form)
in the short-time asymptotic expansion of the heat kernel on $X$. Our method
differs markedly from that of \cite{miat76,miat79-29-249,miat80-260-1}. 
Besides their mathematical interest
these coefficients play an important role in quantum loop effects (such as 
the conformal anomaly), and  in field theory, quantum gravity, and cosmology 
\cite{eliz94,byts96-266-1}.

\section{Acknowledgements}

A.A. Bytsenko wishes to thank CNPq and the Department of Physics of
Londrina University for financial support and kind hospitality. The research
of A.A. Bytsenko was supported in part by Russian Foundation for Basic
Research (grant No. 98-02-18380-a) and by GRACENAS (grant No. 6-18-1997).

\section{Appendix}

The constant $C_G$ in Eq. (2.3) and the Miatello coefficients $a_{2j}$ of the
polynomials $P(r)$ in Eq. (2.5) appear in the statements of Theorems 2.1,
2.2 and 3.1. $C_G$ and $P(r)$ for the various rank 1 simple groups $G$ of
this paper are given in the following table.

\vspace{1.0cm}
{\bf Table}
\begin{flushleft}
\begin{tabular}{||l l l||}
\hline\hline  
${\rm G}$ & ${\rm C}_{{\rm G}}$ & ${\rm P}(r)$ \\
\hline\hline
${\rm SO}_1(n,1),\,\,n\geq 2$ & $\left[2^{2n-4}\Gamma\left(\frac{n}{2}\right)^2
\right]^{-1}$ & 

$\prod_{j=0}^{m-2}\left[r^2+\frac{(2j+1)^2}{4}\right],\,\,
n=2m$\\
& &  
$\prod_{j=0}^{m-1}\left[r^2+j^2\right],\,\,n=2m+1$ \\
${\rm SU}(n,1),\,\,n\geq 2$ & $\left[2^{2n-1}\Gamma(n)^2
\right]^{-1}$ & 
$\prod_{j=1}^{n-1}\left[\frac{r^2}{4}+\frac{(n-2j)^2}{4}
\right]$\\        
${\rm SP}(n,1),\,\,n\geq 2$ & $\left[2^{4n+1}\Gamma(2n)^2
\right]^{-1}$ & 
$\left[\frac{r^2}{4}+\frac{1}{4}\right]\prod_{j=3}^{n+1}
\left[\frac{r^2}{4}+\left(n-j+\frac{3}{2}\right)^2\right]
\left[\frac{r^2}{4}+\left(n-j+\frac{5}{2}\right)^2\right]$
\\        
${\rm F}_{4(-20)}$ & $\left[2^{21}\Gamma(8)^2
\right]^{-1}$ & 
$\left[\frac{r^2}{4}+\frac{1}{4}\right]\left[\frac{r^2}{4}+\frac{9}{4}\right]
\prod_{j=0}^{4}
\left[\frac{r^2}{4}+\left(\frac{2j+1}{2}\right)^2\right]$
\\        
\hline \hline
\end{tabular}
\end{flushleft}

\end{document}